\documentclass[11pt]{amsart}
\newtheorem{theorem}{Theorem}[section]

\theoremstyle{definition}
\newtheorem{remark}{Remark}[section]

\numberwithin{equation}{section}

\newcommand{\Nat}{{\mathbb N}}
\newcommand{\Int}{{\mathbb Z}}

\begin{document}

\title[Divergent correction]{A divergent Vasyunin correction}

\author{Luis B\'{a}ez-Duarte\\ \\12 June 2005}

\date{12 June 2005; from a 31 October 2001 note}



\address{Departamento de Matem\'{a}ticas, 
Instituto Venezolano de Investigaciones
Cient\'{\i}ficas, Apartado 21827, Caracas 1020-A,
 Venezuela}

\email{lbaezd@cantv.net, lbaezd@gmail.net}

\begin{abstract}
V. I. Vasyunin has introduced special sequences of step functions related to the strong Nyman-Beurling criterion that converge pointwise to $1$ in $[1,\infty)$. We show here that the first and simplest such sequence considered by Vasyunin diverges in $L_1((1,\infty),x^{-2}dx)$, which of course precludes the 
$L_2((1,\infty),x^{-2}dx)$-convergence needed for the Riemann hypothesis. Whether all sequences considered by this author also diverge remains an interesting open question.   
\end{abstract}

\maketitle

\section{Introduction}

Any function of the type described by
\begin{equation}\label{natfunction}
\varphi(x)=\sum_{a=1}^N \alpha_a \left[\frac{x}{a}\right], \ \textnormal{with} \ \ \ \sum_{a=1}^N \frac{\alpha_a}{a}=0,
\end{equation}
is called here a \emph{natural function}.

It follows easily from the main theorem in \cite{baez} that the Nyman-Beurling criterion for the Riemann hypothesis in \cite{nyman}, and \cite{beurling} can be strengthened to read as follows: the Riemann hypothesis is equivalent to the existence of a sequence of natural functions converging in 
$L_2((0,\infty),x^{-2}dx)$ to the constant function $1$. 

With this goal in mind V. I. Vasyunin \cite{vasyunin} introduced a very interesting construction of sequences of natural functions, which he called \emph{corrections}, for reasons best understood below in Section \ref{natcorrection} they are \emph{corrections to the natural approximation}. This author first identifies \textit{idempotent} natural functions $f$, call them \emph{seeds}, with the property that $f(x)$ takes only the values $0$ or $1$, then he constructs an appropriate sequence of linear combinations of dilations of $f$, that is equal to $1$ in a large interval $(1,M]$. The classification of seeds in canonical types, which seem to be very scarce, is the object of Vasyunin's paper. 

A new light on the subject is cast by interesting, rather intriguing geometrical research in  A. Borisov's paper \cite{borisov}. This has led the author to resurrect earlier attempts to study the analytical properties of Vasyunin corrections. We show below that the first Vasyunin correction diverges in $L_1((0,1),x^{-2}dx)$, which of course is a lot stronger than divergence in $L_2((0,1),x^{-2}dx)$. It would seem that the proof employed below can be modified and extended to the other Vasyunin corrections, but the difficulties could be considerable. The subject deserves further analytic study: If divergence were generally established we should just have another piece of the Riemann hypothesis puzzle, which in any case does not diminish the intrinsec combinatorial-geometric interest of these constructs. 

\section{Divergence of the first correction} 
Vasyunin first considers the seeds
\begin{equation}\label{firstseed}
f_n(x):=\left[\frac{x}{n}\right]-2\left[\frac{x}{2n}\right],
\end{equation}
and then constructs his \emph{first correction}, the sequence of natural function $\varphi_n$, where each $\varphi_n$ is a sum
\begin{equation}\label{firstcorrection}
\varphi_n=\sum_{k=1}^n c_k f_k
\end{equation}  
with the $c_k$ defined by recurrence as follows
\begin{equation}\label{recurrence} 
c_n=1-\sum_{k=1}^{n-1}c_kf_k(n),
\end{equation}
where the empty sum is $0$ so $c_1=1$. 

It must be remarked that Vasyunin is quite aware of the unsuitable behavior of $\varphi$ as he points out in \cite{vasyunin} that: \textit{However, the described ``correction" of our function near the point $n$ leads to too large perturbations on the interval $(n+1,\infty)$}. This note is written to make precise in some sense these large perturbations.

Vayunin introduced in \cite{vasyunin} other more promising seeds to employ with the same method of construction. Here are Vayunin's second and third seed examples:
$$
g_n(x)=\left[\frac{x}{n}\right]-\left[\frac{x}{n+1}\right]-\left[\frac{x}{n(n+1)}\right],
$$ 
and
$$
h_n(x)=\left[\frac{x}{n}\right]-\left[\frac{x}{n+1}\right]-\left[\frac{x}{2n}\right]+\left[\frac{x}{2(n+1)}\right]-\left[\frac{x}{2n(n+1)}\right].
$$

It seems probable that the proof of the divergence theorem below could be adapted to these cases too, but there is no denying that some important difficulties may appear.

We now state the result announced in the introduction:

\begin{theorem}\label{divthm}
The first correction $\varphi_n$, defined by (\ref{firstcorrection}),
(\ref{recurrence}), and (\ref{firstseed}), diverges in $L_1((0,\infty),x^{-2}dx)$.
\end{theorem}

\begin{proof}
Clearly each $f_k$, and each $\varphi_n$ are right-continuous step functions which are constant in every interval $[a,a+1)$ where $a\in\Nat$. It is clear that
$$
\varphi_1(x)=f_1(x)=[x]-2[x/2]
$$
is $0$ in $[0,1)$ and $1$ in $[1,2)$. Thus 
$$
f_n(x)=0, \ \ \ x\in[0,n),
$$
$$
f_n(x)=1, \ \ \ x\in [n,2n),
$$
It is trivial that all $\varphi_n(x)=0$ for $x\in[0,1)$. We now show by induction that $\varphi_n(x)=1$ for $x\in[1,n+1)$. Clearly this is true for $n=1$; so assume that $\varphi_{n-1}(x)=1$ for $x\in[1,n)$. On account of (\ref{firstcorrection}) $\varphi_n(x)=\varphi_{n-1}(x)+c_n f_n(x)=\varphi_{n-1}(x)=1$ for $x\in[1,n)$, and likewise for $x\in[n,n+1)$ we see from (\ref{firstcorrection}) and (\ref{recurrence}) that
\begin{eqnarray}\nonumber
\varphi_{n}(x)
&=&
\varphi_{n}(n)\\\nonumber
&=&
\varphi_{n-1}(n)+c_n f_n(n)\\\nonumber
&=&
\varphi_{n-1}(n)+c_n\\\nonumber
&=&
\varphi_{n-1}(n)+1-\sum_{k=1}^{n-1}c_k f_k(n)\\\nonumber
&=&
\varphi_{n-1}(n)+1-\varphi_{n-1}(n)=1
\end{eqnarray}
Let now $j_n(m)$ be the jump of $\varphi_n$ at $m\in\Nat$. It is clear that for all $n\in\Nat$ 
\begin{equation}\label{jumps}
j_n(1)=1, \ \ \ j_n(m)=0, \ \ \ 2\leq m\leq n.
\end{equation} 
On the other hand for any $m\in\Nat$ we can write
\begin{eqnarray}\nonumber
j_n(m)
&=&
\sum_{k=1}^n c_k (f_k(m)-f_k(m-0))\\\label{jeq}
&=&
\sum_{1\leq k=\leq n,\ k|m} c_k (-1)^{\frac{m}{k}+1}.
\end{eqnarray}
Denote by $\delta$ the function defined on $\Nat$ given by
$$
\delta(1)=1, \ \ \ \delta(a)=0,\ \ (a>1),
$$
then for any $m\in\Nat$ if one takes $m<n\in\Nat$ the above equation (\ref{jeq}) together with (\ref{jumps}) yield
$$
-\sum_{k|m}c_k (-1)^{\frac{m}{k}}=\delta(m).
$$
This means that the sequence $c_k$ is the Dirichlet inverse of $(-1)^{m+1}$ which obviously exists and can easily be determined because
$$
\sum_{m=1}^\infty \frac{(-1)^{m+1}}{m^s}=(1-2^{1-s})\zeta(s),\ \ \ \
(\Re(s)>0).
$$
For $\Re(s)>1$ the reciprocal of this function can be expressed as 
\begin{eqnarray}\nonumber
\frac{1}{1-\frac{2}{2^s}}\frac{1}{\zeta(s)}
&=&
\sum_{h=0}^\infty\frac{2^h}{(2^h)^s} 
\sum_{j=1}^\infty \frac{\mu(j)}{j^s}\\\nonumber
&=&
\sum_{k=1}^\infty\frac{1}{k^s}\sum_{m|k}\psi(k)\mu\left(\frac{k}{m}\right),
\end{eqnarray}
where
$$
\psi(2^r)=1,\ \ \ (r\in\Int), \ \ \ \textnormal{else}\ \ \ \psi(k)=0,
$$
and clearly
\begin{equation}\label{inv}
c_k=\sum_{j|k}\psi(j)\mu\left(\frac{k}{j}\right).
\end{equation}
Factorize $k=2^r m$, where $r\geq0$ and $m$ is odd, then (\ref{inv}) becomes
\begin{eqnarray}\nonumber
c_k
&=&
\sum_{a=0}^r 2^a \mu(2^{r-a}m)\\\nonumber
&=&
\mu(m)\sum_{a=0}^r  2^a \mu(2^{r-a})\\\nonumber
&=&
\mu(m)\sum_{a=\max(r-1,0)}^r  2^a \mu(2^{r-a})\\\nonumber
&=&
2^{\max(r-1,0)}\mu(m).
\end{eqnarray}
Thus
\begin{equation}\label{coeffexpr}
c_{2^r m}=2^{\max(r-1,0)}\mu(m),\ \ \ (r\geq0, \ m\ \ \textnormal{odd}).
\end{equation}
This shows that the coefficients are too large to even allow $L_1((1,\infty),x^{-2}dx)$-convergence of the $\varphi_n$ since for $n=2^r$, $r\geq1$, we have $c_n=n/2$, so that
\begin{eqnarray}\nonumber
\|\varphi_n-\varphi_{n-1}\|_1 
&=&
\frac{n}{2}\int_1^\infty f_n(x)\frac{dx}{x^2}\\\nonumber
&=&
\frac{n}{2}\int_1^\infty f_1(\frac{x}{n} )\frac{dx}{x^2}\\\nonumber
&=&
\frac{1}{2}\int_{\frac{1}{n}}^\infty f_1(x )\frac{dx}{x^2}\\\nonumber
&=&
\frac{1}{2}\int_1^\infty f_1(x )\frac{dx}{x^2}=C>0.
\end{eqnarray}
\end{proof}
\begin{remark}
From the fact that $\varphi_n(x)-\varphi_{n-1}(x)$ is of one sign it follows from the above that even
$\int_1^\infty(\varphi_n(x)-1)x^{-2}dx\not\rightarrow0$. 
\end{remark}

\section{Connection with the natural approximation}\label{natcorrection}
We should remark in closing that $\varphi_n$ is closely related to the so-called \emph{natural approximation}, namely, the coefficients for $\varphi_n$ in the canonical form (\ref{natfunction}) are $\mu(k)$ at least up to $k=n$, as is indeed the case for any natural function (\ref{natfunction}) with $N\geq n$ with value $1$ in the interval $[1,n+1)$. This is the case also for the corrections arising from more complicated seeds. We now compute all the coefficients of $\varphi_n$ in the canonical form.

\begin{theorem}\label{natapp}
For any $n\in\Nat$
$$
\varphi_n(x)=\sum_{k=1}^n \mu(k)\left[\frac{x}{k}\right]
-2\sum_{\frac{n}{2}< k \leq n}^n c_k\left[\frac{x}{2k}\right].
$$ 
\end{theorem}

\begin{proof}
Using (\ref{firstcorrection}) and (\ref{firstseed}) we have
\begin{eqnarray}\nonumber
\varphi_n(x)
&=&
\sum_{k=1}^n c_k 
\left(\left[\frac{x}{k}\right]-2\left[\frac{x}{2k}\right]\right)\\\nonumber
&=&
\sum_{k\leq n} c_k \left[\frac{x}{k}\right]
-2\sum_{k\leq \frac{n}{2}} c_k \left[\frac{x}{2k}\right]
-2\sum_{\frac{n}{2}< k \leq n}^n c_k\left[\frac{x}{2k}\right].
\end{eqnarray}
We now write
$$
\sum_{k=1}^n\alpha_k  \left[\frac{x}{k}\right] =
\sum_{k\leq n} c_k \left[\frac{x}{k}\right]
-2\sum_{k\leq \frac{n}{2}} c_k \left[\frac{x}{2k}\right]
$$
and try to determine the $\alpha_k$. Use now the expression (\ref{coeffexpr}) for the $c_k$ as follows: If $k$ is odd, then simply $\alpha_k=c_k=\mu(k)$. If $k$ is even, $k=2^r m$, and $r\geq2$, then $\alpha_k=c_k-2c_{\frac{k}{2}}=0=\mu(k)$; whereas for $r=1$ we get $\alpha_k=-\mu(m)=\mu(2m)=\mu(k)$.
\end{proof}

\bibliographystyle{amsplain}

\end{document}